\documentclass{amsart}
\usepackage{graphicx}
\usepackage{amssymb}

\newtheorem{thm}{Theorem}
\newtheorem{lemma}[thm]{Lemma}

\newtheorem{cor}[thm]{Corollary}
\newtheorem{conj}[thm]{Conjecture}
\newtheorem{claim}{Claim}

\newcommand{\Pf}{\noindent \textbf{Proof: }}
\newcommand{\pr}{\noindent \textbf{Proof: }}
\newcommand{\Rmk}{\noindent \textbf{Remark: }}

\newcommand{\R}{{\mathbb R}}

\title{A sufficient condition for intrinsic knotting of bipartite graphs}
\author[S.\ Huck, A.\ Appel, M.-A.\ Manrique, and T.W.\ Mattman]{Sophy Huck, Alexandra Appel, Miguel-Angel Manrique, and Thomas W.\ Mattman}
\address{Department of Mathematics and Statistics,
California State University, Chico,
Chico, CA 95929-0525}
\email{sophy26@hotmail.com}
\email{tinyapple@gmail.com}
\email{TMattman@CSUChico.edu}
\address{Department of Mathematics,
University of Southern California,
3620 South Vermont Ave., KAP 108
Los Angeles, CA 90089-2532}
\email{mmanriqu@usc.edu}
\thanks{The research was supported in part by NSF DMS award 0648764 as part of the Summer 2007 REUT at CSU, Chico}

\dedicatory{Dedicated to the memory of Michael Curtis Wilson.}

\begin{document}

\begin{abstract}
We present evidence in support of a conjecture that a bipartite graph with at least five vertices in each part and $|E(G)| \geq 4 |V(G)| - 17$ is intrinsically knotted. We prove the conjecture for graphs that have exactly
five or exactly six vertices in one part. We also show that there is a constant $C_{n}$ such that a bipartite graph with exactly $n \geq 5$ vertices in one part and  $|E(G)| \geq 4 |V(G)| + C_{n}$ is intrinsically knotted. Finally, we classify bipartite graphs with ten or fewer vertices with respect to intrinsic knotting.
\end{abstract}

\maketitle

\section{Introduction}
We recently discovered~\cite{CMOPRW} that a result of Mader~\cite{M} leads to
a proof of a conjecture of Sachs: A graph on $|V(G)| \geq 7$ vertices with at least $5|V(G)| - 14$ edges is intrinsically knotted. In the current paper, we will give evidence in
support of a similar bound $|E(G)| \geq 4|V(G)| - 17$ that ensures knotting
of bipartite graphs.

Recall that a graph is intrinsically knotted (IK) if every tame embedding of the
graph in $\R^3$ contains a non-trivially knotted cycle.
Since knotless embedding is preserved under edge contraction~\cite{NT}, work of Robertson and Seymour~\cite{RS} shows that this property
is determined by a finite list
of minor minimal IK graphs.
However, determining this list remains difficult.
It is known~\cite{CG,F1,KS,MRS} that $K_7$ and
$K_{3,3,1,1}$ along with any graph obtained from these two by triangle-Y
exchanges is minor minimal with respect to intrinsic knotting.
Foisy~\cite{F2,F3} has shown the existence of several additional
minor minimal IK graphs.

As a complete characterization of IK graphs remains out
of reach for the moment, we propose instead a sufficient condition for
intrinsic knotting that is easy to test. For a graph $G$, let $V(G)$ denote the set of vertices and 
$v(G) = |V(G)|$ the number of vertices.  Similarly, $e(G) = |E(G)|$ will be the number of edges.

\begin{conj} Let $G$ be a bipartite graph with at least five vertices in each part.  If $e(G) \geq 4v(G) - 17$ then $G$ is intrinsically knotted.
\end{conj}

Note that if a bipartite graph has four or fewer vertices in one part,
then it is not IK~\cite{BBFFHL}. Also, any graph can be made bipartite without affecting its topology by adding degree two vertices in the middle of selected edges. So, even for graphs that are not bipartite,
the proposed bound may be more useful
in some instances than the bound $e(G) \geq 5v(G) - 14$ that applies to
all graphs of seven or more vertices.

The conjecture was proved
for graphs with exactly five vertices in one part by~\cite{CHPS}.
We build on their work to show

\begin{thm}
\label{thm56}%
Let $G$ be a bipartite graph with exactly five or exactly six vertices in one part and at least five vertices in the other.
If $e(G) \geq 4v(G) - 17$ then $G$ is IK. Moreover, this is also true if $G$ has exactly seven vertices in each of its two parts.
\end{thm}

For bipartite graphs of ten vertices, the bound $e(G) \geq 4 v(G) - 17
=23$ almost characterizes intrinsic knotting. Let
$K_{n_1,n_2} \setminus m$ denote the set of graphs constructed from the complete bipartite graph $K_{n_1,n_2}$ by removing $m$ edges. We will use $a_1, a_2, \ldots, a_{n_1}$ to denote the vertices in one part (the {\em $a$-vertices}) and $b_1, b_2, \ldots,
b_{n_2}$ those in the other (the {\em $b$-vertices}).

\begin{thm}
\label{thm55}%
A bipartite graph on ten or fewer vertices is not IK unless it has five vertices in each part. A graph $G$ of the form $K_{5,5} \setminus m$  is IK if and only if
\begin{itemize}
\item $G$ has $23$ or more edges, or
\item $G$ is the graph $K_{5,5}$ with the edges
$a_1b_1$, $a_1b_2$, and $a_2b_1$ removed.
\end{itemize}
\end{thm}

We will use the notation $K_{5,5}  \setminus \{a_1b_1, a_1b_2, a_2b_1\}$
to describe the unique IK $K_{5,5} \setminus 3$ graph.

Theorem~\ref{thm55} will be proved in Section 3 below, along with theorems that give sufficient conditions for intrinsic knotting of subgraphs
of $K_{6,6}$ and $K_{7,7}$.
In Section 4, we use these results to prove Theorem~\ref{thm56}
as well as theorems that, together, provide sufficient conditions for intrinsic knotting of subgraphs of any complete bipartite graph.

Although we are unable to prove our conjecture, in Section 5 we will show that a bound of the form $e(G) \geq 4v(G) + C_n$ ensures intrinsic knotting
for bipartite graphs with exactly $n \geq 5$ vertices in one part:

\begin{thm}
\label{thmC}%
Let $a_n$ be defined by the recurrence
$$a_n = \left\lfloor \frac{n(a_{n-1} - 1)}{n-5} \right\rfloor + 1$$
when $n \geq 7$, $a_5 = 5$, and  $a_6 = 7$. Let $C_n = a_n -4n$, for $n \geq 7$, and $C_5 = C_6 = -17$.
Let $G$ be a bipartite graph with exactly $n \geq 5$ vertices in one part
and at least $a_n$ vertices in the other.
If $e(G) \geq 4v(G) + C_n$ then $G$ is IK.
\end{thm}

We remark that our techniques for showing that a graph is IK ultimately come down to showing that it has
one of the graphs obtained from $K_7$ by triangle-Y moves
as a minor.
In particular, we will make use of the graphs $F_9$ and $H_9$ constructed
in this way by Kohara and Suzuki~\cite{KS}.
This means we
have no new examples of minor minimal IK graphs.
In particular, we deduce that there are no new examples to be found
among bipartite graphs on ten or fewer vertices.

This paper is largely an abbreviated version of
\cite{HAM} to which we refer the reader for additional details.

\section{Lemmas}

Our analysis of intrinsic knotting of bipartite graphs is based primarily on the following lemma which follows easily from the Pigeonhole Principle.

\begin{lemma}%
\label{PHPlem}
Let $a$, $b$, $k$, and $m$ be positive integers such that $(k-1)(a + 1) < m + k$. Every graph of the form $K_{a+1,b} \setminus (m+k)$ has a subgraph of
the form $K_{a,b} \setminus m$.
\end{lemma}

Recall that any graph obtained from the graph $G$ by a sequence of
edge deletions or contractions is called a minor of $G$. It follows from
\cite{NT} that if $H$ is IK and $H$ is a minor of $G$, then $G$ is
also IK. Since subgraphs are examples of minors, Lemma~\ref{PHPlem} implies 

\begin{lemma}%
\label{lemmain}
Let $a$, $b$, $k$, and $m$ be positive integers such that $(k-1)(a + 1) < m + k$. If every graph of the form $K_{a,b} \setminus m$ is IK, then the same is true of every graph of
the form $K_{a+1,b} \setminus (m+k)$.
\end{lemma}

We also make note of a useful lemma due to \cite{BBFFHL,OT}.
Let $K_2 + G$ denote the join of $G$ and $K_2$, the complete graph on
$2$ vertices.

\begin{lemma}%
\label{Flemlem}
$K_2 + G$ is IK if and only if $G$ is non-planar.
\end{lemma}

\section{Subgraphs of $K_{5,5}$, $K_{6,6}$, and $K_{7,7}$}

In this section we prove Theorem~\ref{thm55}, characterizing intrinsic knotting of bipartite graphs on 10 or fewer vertices, as well as
theorems that give sufficient conditions for a subgraph of $K_{6,6}$ or $K_{7,7}$ to be IK.

\setcounter{thm}{2}

\begin{thm}
A bipartite graph on ten or fewer vertices is not IK unless it has five vertices in each part. A graph $G$ of the form $K_{5,5} \setminus m$  is IK if and only if
\begin{itemize}
\item $G$ has $23$ or more edges, or
\item $G = K_{5,5} \setminus \{a_1b_1, a_1b_2, a_2b_1\}$.
\end{itemize}
\end{thm}

\setcounter{thm}{7}

\pr
It was shown in \cite{BBFFHL} that a bipartite graph with $4$ or fewer vertices in one part is not intrinsically knotted. In \cite{CMOPRW} it is argued that all graphs of the form $K_{5,5} \setminus 2$ are IK, so any graph $K_{5,5} \setminus m$ with $23$ or more edges is IK.

There are
four graphs of the form $K_{5,5} \setminus 3$. Three of these are subgraphs  of a graph of the form $K_2+H$ where $H$ is a planar graph and are not IK by Lemma~\ref{Flemlem}. The fourth, $K_{5,5} \setminus \{a_1b_1, a_1b_2, a_2b_1\}$
was shown in Figure 8 of \cite{MOR}  to have the IK graph $H_9$
as a minor and is therefore also intrinsically knotted. Here, $H_9$ is
one of the graphs obtained by triangle-Y substitution on $K_7$; see
\cite{KS}.

All but one of the $K_{5,5} \setminus 4$ graphs are subgraphs
of an unknotted $K_{5,5} \setminus 3$ and are, therefore, also not IK.
The remaining graph, $K_{5,5} \setminus \{a_1b_1, a_1b_2, a_2b_1, a_2b_2 \}$ is, again, a subgraph of a $K_2+H$ graph with $H$ planar, so not IK.
Since no $K_{5,5} \setminus 4$ graph is IK, no
$K_{5,5} \setminus m$ graph with $m \geq 4$ is IK.
\qed

\bigskip

\begin{figure}[ht]
\begin{center}
\includegraphics[scale=0.7]{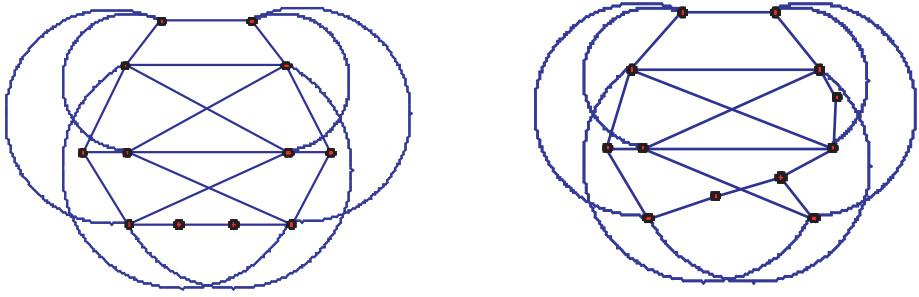}
\caption{The $K_{6,6}\setminus 12$ graphs $H_{66}$ (at left) and $F_{66}$ (at right).}\label{fig1}
\end{center}
\end{figure}

The following two theorems make use of graphs $F_{66}$ and $H_{66}$ (see Figure~\ref{fig1}) of the form $K_{6,6} \setminus 12$ that have as minors the
graphs $F_9$ and $H_9$ obtained from $K_7$ by triangle-Y moves in \cite{KS}.  Since $F_9$ and $H_9$ are IK, their expansions $F_{66}$ and $H_{66}$ are also.

\begin{thm} \label{thm665}%
Every graph of the form $K_{6,6} \setminus 5$ is IK.
\end{thm}

\pr
We will list all ways to remove five edges from the complete bipartite graph $K_{6,6}$ and demonstrate that every case but one (the first case below) has $H_{66}$ as a subgraph. However, the one exceptional case instead has $F_{66}$ as a minor. Thus, in any case, a $K_{6,6} \setminus 5$ graph has an IK minor and is, therefore, IK.

To determine all possible ways to remove five edges from $K_{6,6}$, let $a_1, \ldots , a_6$ be the vertices in one part and $b_1, \ldots , b_6$ those in the other part. Now consider any partition of 5 and take the $i$-th element in a partition to be the number of edges removed from the $i$-th vertex in the first part of $K_{6,6}$.  Likewise, consider another partition and allow its entries to correspond to the number of edges removed from the other part of $K_{6,6}$.  Some combinations cannot be realized as  $K_{6,6}  \setminus 5$ graphs (for instance, the pairing $\{ \{5 \}, \{2,2,1\}\}$ cannot be constructed).  Observe that there may be more than one graph for a given pairing.  Below we indicate the pair of partitions of 5 for each of the twenty
$K_{6,6} \setminus 5$ graphs.

Many of these graphs can also be shown to be IK using Corollary 2.5 of \cite{CHPS} which states that a $K_{6,5}$ graph with two or fewer edges removed is IK. In those cases, the list below specifies which vertex to remove to arrive at such a subgraph of $K_{6,5}$.  

\begin{enumerate}
\item
$ K_{6,6} \setminus \{a_1b_1,a_2b_2,a_3b_3,a_4b_4,a_5b_5\}$;    $\{ \{ 1,1,1,1,1\} ,\{ 1,1,1,1,1\} \}$
\item $  K_{6,6} \setminus \{ a_1b_5, a_2b_5, a_4b_4, a_5b_1, a_6b_6 \}$;    $\{ \{ 2,1,1,1\} , \{ 1,1,1,1,1\} \}$
\item $ K_{6,6} \setminus \{ a_1b_5, a_2b_5, a_4b_1, a_5b_1, a_6b_6 \}$;   $\{ \{ 2,2,1 \} , \{ 1,1,1,1,1\} \}$
\item $K_{6,6} \setminus \{a_1b_5, a_2b_5, a_5b_1, a_5b_2, a_6b_6\}$ ;    $\{ \{2,1,1,1\} , \{ 1,1,2,1\}  \}$ 
\item $K_{6,6} \setminus \{a_2b_5, a_4b_1, a_5b_1, a_5b_2, a_6b_6\}$ ;    $\{ \{2,1,1,1\} , \{ 1,2,1,1\}  \}$ 
 \item $K_{6,6} \setminus \{a_1b_5, a_2b_5, a_5b_1, a_5b_2, a_6b_6\}$ ;    $\{ \{2,1,2\} , \{ 1,2,1,1\}  \}$ 
 \item $K_{6,6} \setminus \{a_1b_5, a_2b_5, a_5b_1, a_5b_2, a_5b_3\}$ ;    $\{ \{2,1,1,1\} , \{ 1,1,3\}  \}$ 
Remove vertex $a_3$ to obtain a $K_{6,5} \setminus 2$ graph which we know to be IK \cite[Corollary 2.5]{ CHPS}.
\item $K_{6,6} \setminus \{a_1b_4, a_1b_5, a_5b_1, a_5b_2, a_5b_3\}$ ;    $\{ \{2,3\} , \{ 1,1,1,1,1\}  \}$ 
Remove vertex $a_3$ \cite{CHPS}.
 \item $K_{6,6} \setminus \{a_1b_5, a_4b_1, a_4b_4, a_4b_5, a_5b_1\}$ ;    $\{ \{3,1,1\} , \{ 1,2,2\}  \}$ 
Remove vertex $a_1$ \cite{CHPS}.
\item $K_{6,6} \setminus \{a_4b_1, a_4b_4, a_4b_5, a_5b_1, a_6b_6\}$ ;    $\{ \{3,1,1\} , \{ 1,1,2,1\}  \}$ 
Remove vertex $a_1$ \cite{CHPS}.
\item $K_{6,6} \setminus \{a_4b_4, a_5b_1, a_5b_2, a_5b_3, a_6b_6\}$ ;    $\{ \{3,1,1\} , \{ 1,1,1,1,1\}  \}$ 
Remove vertex $a_1$ \cite{CHPS}.
\item $K_{6,6} \setminus \{a_1b_4, a_1b_5, a_5b_1, a_5b_2, a_5b_4\}$ ;    $\{ \{2,3\} , \{ 1,2,1,1\}  \}$ 
Remove vertex $a_3$ \cite{CHPS}.
\item $K_{6,6} \setminus \{a_1b_4, a_5b_1, a_5b_2, a_5b_3, a_5b_4\}$ ;    $\{ \{1,4\} , \{ 2,1,1,1\}  \}$ 
Remove vertex $a_2$ \cite{CHPS}.
\item $K_{6,6} \setminus \{a_5b_1, a_5b_2, a_5b_3, a_5b_4, a_6b_6\}$ ;    $\{ \{4,1\} , \{ 1,1,1,1,1\}  \}$ 
Remove vertex $a_1$ \cite{CHPS}.
\item $K_{6,6} \setminus \{a_1b_1, a_1b_2, a_1b_3, a_1b_4, a_1b_5 \}$ ;    $\{ \{ 5\} , \{ 1,1,1,1,1\}  \}$ 
Remove the vertex $a_1$ from the original graph to obtain a complete $K_{6,5}$ graph.  This graph has a $K_{5,5}$ minor, which is known to be IK \cite{S}.
\item $K_{6,6} \setminus \{a_1b_4, a_4b_1, a_4b_5, a_5b_1, a_5b_4\}$ ;    $\{ \{1,2,2\} , \{ 2,2,1\}  \}$  \item $K_{6,6} \setminus \{a_1b_4, a_4b_1, a_5b_1, a_5b_4, a_6b_6\}$ ;    $\{ \{1,2,1,1\} , \{ 2,2,1\}  \}$ 
\item $K_{6,6} \setminus \{a_1b_4, a_1b_5, a_2b_5, a_4b_4, a_4b_5\}$ ;    $\{ \{2,2,1\} , \{ 2,3\}  \}$ 
Remove vertex $b_2$ \cite{CHPS}.
\item $K_{6,6} \setminus \{a_1b_4, a_1b_5, a_4b_4, a_4b_5, a_6b_6\}$ ;    $\{ \{2,2,1\} , \{ 2,2,1\}  \}$ 
\item $K_{6,6} \setminus \{a_1b_4, a_4b_1, a_4b_4, a_4b_5, a_5b_4\}$ ;    $\{ \{1,3,1\} , \{ 1,3,1\}  \}$ 
Remove vertex $a_2$ \cite{CHPS}. 
\end{enumerate}

 \qed

\begin{thm} \label{thm666}%
A graph of the form $K_{6,6} \setminus 6$ is IK provided it is not the graph $G_{666} = K_{6,6} \setminus \{a_1b_1, a_2b_2, a_3b_3, a_4b_4, a_5b_5, a_6b_6 \}$.
\end{thm}

\Rmk  We do not know whether $G_{666}$ is IK or not. Using a computer algebra system, we have verified that none of the graphs obtained by triangle-Y exchanges on $K_7$ or $K_{3,3,1,1}$ is a minor of $G_{666}$. So, if $G_{666}$ is IK, it or one of its minors would be a new example of
a minor minimal IK graph.

\bigskip

\pr
Let $G$ be a graph of the form $K_{6,6} \setminus 6$. If $G$ has a vertex
from which three edges have been removed (compared to the complete bipartite graph $K_{6,6}$), we may delete that vertex to obtain an IK subgraph of the form
$K_{6,5} \setminus m$ with $m \leq 3$ \cite{CHPS}.
So, we may assume that each vertex in $G$ has at most two edges removed.
Thus, it will be enough to consider partitions of six into parts of size one or two. Below we 
list the seventeen graphs of this form.

\begin{enumerate}
\item $K_{6,6} \setminus \{a_1b_1, a_2b_2, a_3b_3, a_4b_4, a_5b_5, a_6b_6\}$ ;    $\{ \{1,1,1,1,1,1\} , \{ 1,1,1,1,1,1\}  \}$ 
\item $K_{6,6} \setminus \{a_1b_4, a_2b_5, a_4b_1, a_4b_4, a_5b_2, a_6b_6\}$ ;    $\{ \{2,1,1,1,1\} , \{ 1,1,1,1,1,1\}  \}$ 
\item $K_{6,6} \setminus \{a_1b_4, a_2b_5, a_3b_5, a_4b_1, a_5b_4, a_6b_6\}$ ;    $\{ \{2,2,1,1\} , \{ 1,1,1,1,1,1\}  \}$ 
\item $K_{6,6} \setminus \{a_1b_5, a_2b_5, a_4b_1, a_4b_4, a_5b_3, a_6b_6\}$ ;    $\{ \{2,1,1,1,1\} , \{ 1,1,2,1,1\}  \}$ 
\item $K_{6,6} \setminus \{a_1b_4, a_2b_5, a_4b_1, a_4b_4, a_5b_2, a_6b_6\}$ ;    $\{ \{2,1,1,1,1\} , \{ 1,2,1,1,1\}  \}$ 
\item $K_{6,6} \setminus \{a_1b_4, a_3b_5, a_4b_1, a_5b_1, a_5b_4, a_6b_6\}$ ;    $\{ \{2,2,1,1\} , \{ 1,2,1,1,1\}  \}$ 
\item $K_{6,6} \setminus \{a_1b_4, a_2b_5, a_3b_5, a_5b_1, a_5b_4, a_6b_6\}$ ;    $\{ \{2,1,2,1\} , \{ 1,2,1,1,1\}  \}$ 
\item $K_{6,6} \setminus \{a_1b_4, a_2b_5, a_3b_5, a_4b_1, a_5b_1, a_5b_4\}$ ;    $\{ \{2,2,2\} , \{ 1,2,1,1,1\}  \}$ 
\item $K_{6,6} \setminus \{a_2b_5, a_3b_5, a_4b_1, a_4b_4, a_5b_1, a_5b_2\}$ ;    $\{ \{2,2,1,1\} , \{ 1,2,1,2\}  \}$ 
\item $K_{6,6} \setminus \{a_3b_1, a_3b_5, a_5b_2, a_5b_4, a_6b_3, a_6b_6\}$ ;    $\{ \{2,2,2\} , \{ 1,1,1,1,1,1\}  \}$ 
\item $K_{6,6} \setminus \{a_1b_4, a_2b_5, a_3b_5, a_4b_4, a_5b_1, a_5b_2\}$ ;    $\{ \{2,1,1,2\} , \{ 1,1,2,1,1\}  \}$ 
\item $K_{6,6} \setminus \{a_1b_4, a_2b_5, a_4b_1, a_4b_5, a_5b_2, a_5b_4\}$ ;    $\{ \{2,1,2,1\} , \{ 1,2,1,2\}  \}$ 
\item $K_{6,6} \setminus \{a_1b_4, a_1b_5, a_2b_5, a_5b_1, a_5b_4, a_6b_6\}$ ;    $\{ \{2,2,1,1\} , \{ 1,2,2,1\}  \}$ 
\item $K_{6,6} \setminus \{a_1b_4, a_2b_5, a_4b_1, a_4b_5, a_5b_1, a_5b_4\}$ ;    $\{ \{2,2,2\} , \{ 1,2,2,1\}  \}$ 
\item $K_{6,6} \setminus \{a_1b_5, a_1b_4, a_4b_4, a_4b_5, a_5b_1, a_6b_6\}$ ;    $\{ \{2,2,1,1\} , \{ 2,2,1,1\}  \}$ 
\item $K_{6,6} \setminus \{a_1b_5, a_2b_5, a_4b_1, a_4b_4, a_5b_1, a_5b_4\}$ ;    $\{ \{2,2,2\} , \{ 2,2,1,1\}  \}$ 
\item $K_{6,6} \setminus \{a_1b_4, a_1b_5, a_4b_1, a_4b_5, a_5b_1, a_5b_4\}$ ;    $\{ \{2,2,2\} , \{ 2,2,2\}  \}$ 
\end{enumerate}

All but two of the graphs have $H_{66}$ as a subgraph. The two exceptions are
the eleventh graph  $K_{6,6} \setminus \{a_1b_4,a_2b_5,a_3b_5,a_4b_4,a_5b_1,a_5b_2 \}$, which instead has $F_{66}$ as a subgraph, and the first graph $G_{666}$.
\qed

\begin{thm} \label{thm7710}
Any graph of the form $K_{7,7} \setminus 10$ is IK.
\end{thm}

\pr
Let $G$ be a graph of the form $K_{7,7} \setminus 10$. If $G$ has a vertex from which three edges have been removed, that vertex may be deleted to obtain a $K_{7,6} \setminus m$ subgraph with $m \leq 7$. We will show that such graphs are IK in Theorem~\ref{thmK6}, below.  So, we may assume that each vertex of $G$ has at most two edges removed.

Suppose there are vertices $a$ and $b$, one in each part, each having two edges removed, and such that $ab$ is an edge of $G$.
We will argue that $G$ is IK in this case. This will complete the argument as there must be two vertices with these properties. Indeed, as there are 10 edges removed and $G$ has 7 vertices on each side, there must be at least three vertices on each side with two edges removed. Call them $a_1$, $a_2$, $a_3$, $b_1$, $b_2$, and $b_3$. If there are no pair $a_i$, $b_j$ such that $a_ib_j$ is in $G$, then each of these six vertices in fact has three edges removed, violating an earlier assumption.

So, we may assume that there are vertices $a_1$ and $b_1$, each with two
edges removed and such that $a_1b_1$ is an edge of $G$.
Then removing
$a_1$ and $b_1$ from $G$ results in a $K_{6,6} \setminus 6$ graph. If that
graph is not $G_{666}$, then, by Theorem~\ref{thmK6}, $G$ has an IK subgraph and is itself IK.

Thus, we may assume that removing $a_1$ and $b_1$ results in $G_{666}$.
Label the remaining vertices of $G$ so that the edges $a_2b_2$, $a_3b_3$, ..., $a_7b_7$ are not in $G$. Let $a_1b_2$ and $a_1b_3$ denote the edges removed from $a_1$. That is, $a_1b_2$ and $a_1b_3$ are not edges of $G$.
Then, up to relabeling, there are three possibilities for the edges removed from $b_1$.
In each case we will show how removing a different pair of vertices
from $G$ will result in a $K_{6,6} \setminus 6$ subgraph other than $G_{666}$.

The first case is that $a_2b_1$ and $a_3b_1$ are the edges removed from $b_1$. Then, we
can obtain a $K_{6,6} \setminus 6$ graph different from $G_{666}$ by removing the vertices $a_4$ and $b_2$ instead of $a_1$ and $b_1$. That means, in this case, $G$ has an IK subgraph and is also IK.

Next, suppose
that $a_2b_1$ and $a_4b_1$ are the edges removed from $b_1$. Here, we can
obtain a $K_{6,6} \setminus 6$ graph other than $G_{666}$ by removing
instead the vertices $a_1$ and $b_3$.

Finally, if $a_4b_1$ and $a_5b_1$ are the edges removed from $b_1$, then we should instead remove the vertices $a_4$ and $b_2$. \qed

\section{Sufficient conditions for subgraphs of any complete bipartite graph}

In this section we present some general results that, taken together with those of the previous section, provide sufficient conditions for intrinsic knotting of a subgraph of $K_{a,b}$ for any $a,b \geq 6$. Together with \cite{CHPS} (who treat the case where one of $a$ or $b$ is exactly five) and \cite{BBFFHL} (who show that a bipartite graph with at most four vertices in one part is not IK) this gives a collection of sufficient conditions for subgraphs of any complete bipartite graph. We will
also prove Theorem~\ref{thm56} which shows that our conjecture holds for
subgraphs of $K_{7,7}$ as well as bipartite graphs having exactly five or six vertices in one part.

\begin{thm} \label{thmK6}%
All graphs of the form
$K_{6+n, 6} \setminus (2n + 5)$ with $n \geq 1$ are IK.
\end{thm}

\Pf: We use induction on $n$.
For the base step, let $n = 1$. We will look at four cases.
In each case we will see that there is a way to remove a vertex to
obtain an IK $K_{6,6} \setminus m$ subgraph.

	Case 1:  Consider the $K_{7, 6} \setminus 7$ graph where each of the seven $a$-vertices has exactly one edge removed, and exactly one of the six $b$-vertices (call it $b_1$) has exactly 2 edges removed.

We can remove one vertex to obtain a $K_{6, 6}\setminus 6$ subgraph, all of which have been shown (Theorem~\ref{thm666}) to be IK, except  $G_{666}$. We can avoid $G_{666}$  by not removing one of the two $a$-vertices that is connected to $b_1$ in the complement graph.  Removing any of the other $a$-vertices will leave us with a $K_{6, 6} \setminus 6$ graph that is IK.

	Case 2:  Consider the $K_{7, 6} \setminus 7$ graphs where each of the seven $a$-vertices has exactly one edge removed and at least one of the $b$-vertices has three or more edges removed.  Removing any $a$-vertex results in a $K_{6, 6}\setminus 6$ subgraph that is known to be IK by Theorem~\ref{thm666}.

	Case 3:  Consider the $K_{7, 6} \setminus 7$ graphs where each of the seven $a$-vertices has exactly one edge removed and at least two $b$-vertices have at least two missing edges.  Removing any $a$-vertex results in a $K_{6, 6}\setminus 6$ subgraph that is known to be IK.

	Case 4:  Consider the $K_{7, 6} \setminus 7$ graphs where one or more of the seven $a$-vertices has more than one edge removed.  By removing one of those vertices, we will be left with a $K_{6, 6} \setminus m$ subgraph where $m \leq 5$, all of which are known to be IK by Theorem~\ref{thm665}.

As these four cases cover all possibilities, all $K_{7, 6} \setminus 7$ graphs are IK.

Now, for the inductive step, let $n \geq 1$ and assume all $K_{6+n, 6} \setminus (2n + 5)$ graphs are IK. Using Lemma~\ref{lemmain}
we see that every $K_{6+n+1, 6} \setminus (2(n+1) + 5)$ graph G is also IK.

    Thus, for $n \geq 1$, every $K_{6+n, 6} \setminus (2n + 5)$ graph is IK. \qed

\bigskip

Theorem~\ref{thm665} shows that the above result is also true when $n = 0$.

\begin{cor} Let $G$ be a bipartite graph $G$ with exactly six vertices in one part and at least six vertices in the other part.
If $e(G) \geq 4v(G) - 17$, then $G$ is IK.
\end{cor}

We can now prove Theorem~\ref{thm56}.

\setcounter{thm}{1}

\begin{thm}
Let $G$ be a bipartite graph with exactly five or exactly six vertices in one part and at least five vertices in the other.
If $e(G) \geq 4v(G) - 17$ then $G$ is IK. Moreover, this is also true if $G$ has exactly seven vertices in each of its two parts.
\end{thm}

\setcounter{thm}{12}

\Pf
The case of exactly five vertices was proved in \cite{CHPS}. The case of exactly six vertices is the corollary above while the subgraph of $K_{7,7}$ case is Theorem~\ref{thm7710}. \qed

\begin{thm} \label{thmK7}%
All graphs of the form
$K_{7+n, 7} \setminus (2n + 10)$ with $n \geq 1$ are IK.
\end{thm}
\Pf     We use induction on $n$.

For the base case, let $n = 1$. This gives us a $K_{8, 7} \setminus 12$ graph. From Lemma~\ref{lemmain} and Theorem~\ref{thm7710} we know that all $K_{8, 7} \setminus 12$ graphs are IK.

Now, for the inductive step, let $n \geq 1$ and assume all $K_{7+n, 7} \setminus (2n + 10)$ graphs are IK. By Lemma~\ref{lemmain}, every $K_{7+n+1, 7} \setminus (2(n + 1) + 10)$ graph is also IK. \qed

\bigskip
Note that Theorem~\ref{thm7710} corresponds to the case $n=0$ of this theorem.

Although we cannot verify our conjecture in this case, Theorem~\ref{thmK7} does allow us to improve on the previous best known bound (for graphs in general) \cite{CMOPRW} of $e(G) \geq 5v(G) - 14$:

\begin{cor}
Let $G$ be a bipartite graph with exactly seven vertices in one part and at least seven vertices in the other part. If $e(G) \geq 5v(G) - 31$ then $G$ is IK.
\end{cor}

\begin{thm} \label{thmK8}
All graphs of the form $K_{8+n, 8} \setminus (2n + 15)$ with $n \geq 1$ are IK.
\end{thm}
\Pf     We use induction on $n$.

For the base case, let $n = 1$. This gives us a $K_{9, 8} \setminus 17$ graph $G$. Applying Lemma~\ref{lemmain}, we obtain a $K_{9, 7} \setminus 14$ subgraph, which is IK by Theorem~\ref{thmK7}. Therefore $G$ is IK.

Now, for the inductive step, let $n \geq 1$ and assume all $K_{8+n, 8} \setminus (2n + 15)$ graphs are IK. From Lemma~\ref{lemmain}, we see that every $K_{8+n+1, 8} \setminus (2(n + 1) + 15)$ graph is also IK. \qed

\bigskip

Although we cannot prove the analogue of Theorem~\ref{thmK8} when $n = 0$,  it follows from Theorem~\ref{thmK7} with $n=1$ and Lemma~\ref{lemmain} with $k=2$ that all $K_{8,8} \setminus 14$ graphs are IK.

\begin{thm} \label{thmKaa}
Any graph of the form $K_{a, a} \setminus (6a -34)$, where $a \geq 9$, is IK.
\end{thm}

\Pf   We use induction on $a$.

For the base case, let $a = 9$. Consider a $K_{9, 9} \setminus 20$ graph. By Theorem~\ref{thmK8}, every $K_{9,8} \setminus 17$ graph is IK. Then, by Lemma~\ref{lemmain}, every $K_{9, 9} \setminus 20$ is also IK.
	
Now, for the inductive step, let $a \geq 9$ and assume all $K_{a, a} \setminus (6a -34)$ graphs are IK. We will show that every $K_{a+1, a+1} \setminus (6(a + 1) - 34)$ graph is IK.

Applying Lemma~\ref{lemmain} shows that every $K_{a+1,a} \setminus (6a - 31)$ graph is IK.  Another application shows that all $K_{a+1,a+1} \setminus (6a - 28)$ graphs are IK.  Since $6a - 28 = 6(a+1) - 34$, every $K_{a+1, a+1} \setminus (6(a+1) - 34)$ graph is IK.

Thus, by induction every $K_{a, a} \setminus (6a -34)$ graph is IK for all $a \geq 9$. \qed

\begin{thm}
Any graph of the form $K_{a+n, a} \setminus  (3n + 6a-34)$, where
 $a \geq 9$ and $n \geq 0$, is IK.
\end{thm}

\Pf     We use induction on $n$.

For the base case, let $n = 0$.  This gives us a $K_{a,a} \setminus (6a-34)$ graph,  which is IK by Theorem~\ref{thmKaa}.

Now, for the inductive step, assume that $K_{a+n, a} \setminus (3n + 6a-34)$ is IK.  Let's look at $K_{a+n+1,a} \setminus (3(n + 1) + 6a-34)$, or equivalently, $K_{a+n+1,a} \setminus (3n + 6a-31)$.  Using Lemma~\ref{lemmain}, $K_{a+n+1,a} \setminus (3n + 6a-31)$ is IK provided $3n + 6a-31 > 2(a + n + 1)$.  This holds when $a = 9$ since $3n + 23 > 2(n + 10)$.  Each time $a$ increases by $1$, $3n + 6a-31$ increases by $6$ while $2(a + n + 1)$ increases by $2$.  By induction, the inequality in Lemma~\ref{lemmain} holds for all $a$ and therefore each $K_{a+n+1,a} \setminus (3n+6a-31)$ graph is IK. \qed

\section{Proof of Theorem~\ref{thmC}}

\setcounter{thm}{3}

In this section we prove
\begin{thm}
Let $a_n$ be defined by the recurrence
$$a_n = \left\lfloor \frac{n(a_{n-1} - 1)}{n-5} \right\rfloor + 1$$
when $n \geq 7$, $a_5 = 5$, and  $a_6 = 7$. Let $C_n = a_n -4n$, for $n \geq 7$, and $C_5 = C_6 = -17$. Let $G$ be a bipartite graph with exactly $n \geq 5$ vertices in one part and at least $a_n$ vertices in the other. If $e(G) \geq 4v(G) + C_n$ then $G$ is IK.
\end{thm}

\Pf
Since $C_5 = C_6 = -17$ was proved in Theorem~\ref{thm56}, we will assume that $n \geq 7$. It is enough to show that

\begin{claim}All $K_{n,a} \setminus ((n-4)a -a_n)$ graphs with $a \geq a_n$ are IK.
\end{claim}

Let us first observe that this will complete the argument. Indeed, let $G$ be a $K_{n,a} \setminus m$ graph where $a \geq a_n$ and suppose $e(G) \geq 4v(G) + C_n$. That is, $na-m \geq 4(n+a) + C_n$. Then,
\begin{eqnarray*}
m & \leq & na - 4(n+a) - C_n\\
  & = & na - 4(n+a) - (a_n-4n) \\
  & = & (n-4)a -a_n.
\end{eqnarray*}
So, if we can establish our claim, it will follow that $G$ is IK, as required.

We prove the claim by induction. The case where $n = 6$ corresponds to Theorem~\ref{thmK6}. So, let $n \geq 7$.

We first observe that the claim is valid when $a = a_n$. Indeed, by the inductive hypothesis, all $K_{n-1,a_n} \setminus(n-1-4)a_n - a_{n-1}$ graphs are IK. Then, applying Lemma~\ref{lemmain}, it follows that all $K_{n,a_n} \setminus (n-5)a_n$ graphs are IK, as required.

Now use Lemma~\ref{lemmain} to show, by induction, that the claim holds for all $a \geq a_n$. \qed

\end{document}